\documentclass[12pt]{amsart}

\newcommand \bpig {{\textup B}_{\pi}(G)}
\newcommand \bpi {{\textup B}_{\pi}}

\newcommand \ipig {{\textup I}_{\pi}(G)}
\newcommand \npig {{\textup N}_{\pi}(G)}
\newcommand \npi {{\textup N}_{\pi}}
\newcommand \irr {\textup{Irr}}
\newcommand \irrg {\textup{Irr}(G)}
\newcommand \sbs {\subseteq}
\newcommand \dvd {\hbox {\Big|}}
\newcommand \ndvd {\hbox {\big/}\kern-5pt\dvd}
\newcommand \nrml {\lhd}
\def \< {\langle}
\def \> {\rangle}
\newcounter{palistlocal}

\title[Two distinct canonical sets of lifts of Brauer characters]{A
  construction of two distinct canonical sets of lifts of Brauer
  characters of a $p$-solvable group}

\keywords{Brauer characters, lifts, p-solvable} \subjclass{Primary 20C20}

\author{James P. Cossey}
\address{Department of Mathematics, University of Arizona, Tucson,
AZ, 85721} \email{cossey@math.arizona.edu}

\begin{document}

\begin{abstract}  In \cite{npiconst},
Navarro defines the set $\mbox{Irr}(G \mid Q, \delta) \subseteq \mbox{Irr}(G)$, where $Q$
is a $p$-subgroup of a $p$-solvable group $G$, and shows that if
$\delta$ is the trivial character of $Q$, then $\mbox{Irr}(G \mid Q, \delta)$ provides a set of canonical lifts
of ${\textup{IBr}}_p(G)$, the irreducible Brauer characters with vertex $Q$.  Previously, in \cite{bpi}, Isaacs defined a
canonical set of lifts $\bpig$ of $\ipig$.  Both of these results extend the Fong-Swan Theorem to $\pi$-separable groups, and both construct canonical sets of lifts of the generalized Brauer characters.  It is known that
in the case that $2 \in \pi$, or if $|G| $ is odd, we
have $\bpig = \mbox{Irr}(G \mid Q, 1_Q)$. In this note we give a
counterexample to show that this is not the case when $2 \not\in \pi$.
It is known that if $N \nrml G$ and $\chi \in \bpig$, then the
constituents of $\chi_N$ are  in $\bpi(N)$.  However, we use the same counterexample to show that if $N \nrml G$, and $\chi
\in \mbox{Irr}(G\mid Q, 1_Q)$ is such that $\theta \in \mbox{Irr}(N)$ and $[\theta, \chi_N]
\neq 0$, then it is not necessarily the case that $\theta \in \textup{Irr}(N)$ inherits this property.

\end{abstract}

\maketitle

{\bf{1. Introduction.}} \setcounter{equation}{0}  Let $p$ be a prime, $G$ a finite $p$-solvable group, and for a class function $\alpha$ of $G$, let $\alpha^0$ denote the restriction of $\alpha$ to the $p$-regular elements of $G$.  The celebrated Fong-Swan Theorem asserts that if $\varphi$ is an irreducible Brauer character of $G$ for the prime $p$, then there necessarily exists an ordinary irreducible character $\chi$ of $G$ such that $\chi^0 = \varphi$.  Such a character $\chi$ is called a lift of $\varphi$.  In \cite{bpi}, Isaacs constructs a canonical set of lifts ${\textup B}_{p'}(G) \subseteq \textup{Irr}(G)$, such that for each Brauer character $\varphi$ of $G$, there is a unique lift of $\varphi$ in ${\textup B}_{p'}(G)$.  Similarly, in \cite{npiconst}, Navarro constructs another canonical set $\textup{Irr}(G | Q, 1_Q)$ (which we will denote by ${\textup N}_{p'}(G)$) of lifts of the irreducible Brauer characters of $G$.  Navarro conjectured but did not prove that ${\textup B}_{p'}(G) \neq {\textup N}_{p'}(G)$.  In this paper we give conditions under which ${\textup B}_{p'}(G) = {\textup N}_{p'}(G)$, and give an example to show that these sets need not be equal in general.

Let $\pi$ be a set of primes, and denote by $\pi'$ the complement of $\pi$.  (In the classical case, $\pi$ is the complement of the prime $p$.)  Let $G$ be a $\pi$-separable group.  In
\cite{gajen}, Gajendragadkar constructs a certain class of characters, called the $\pi$-special
characters of $G$.  An irreducible character $\chi$ is $\pi$-special if (a)
$\chi(1)$ is a $\pi$-number and (b) if whenever $S \nrml \nrml G$ and $\theta
\in \mbox{Irr}(S)$ lies under $\chi$, then the order of $\mbox{det}(\theta)$, as a
character of $S/S'$, is a $\pi$-number.  These $\pi$-special
characters are known to have many interesting properties, and they are
necessary for constructing the sets of characters under discussion in
this paper.  In particular, if $\alpha \in \mbox{Irr}(G)$ is $\pi$-special
and $\beta \in \mbox{Irr}(G)$ is $\pi'$-special, then $\alpha \beta \in
\mbox{Irr}(G)$ and this factorization is unique.  If $\chi \in
\mbox{Irr}(G)$ can be written in this factored form, we say that
$\chi$ is $\pi$-factorable and we let $\chi_{\pi}$ and $\chi_{\pi'}$
be the $\pi$-special and $\pi'$-special factors of $\chi$, respectively.  In addition, if $M \triangleleft G$ has $\pi'$-index in $G$, and if $\gamma \in
\mbox{Irr}(M)$ is $\pi$-special and invariant in $G$, then there is a unique $\pi$-special
character $\alpha
\in \mbox{Irr}(G)$ such that $\alpha$ extends $\gamma$.  If $\gamma
\in \mbox{Irr}(M)$ is $\pi$-special and $M \nrml G$ with $G/M$ is a $\pi$-
group, then every character $\alpha \in \mbox{Irr}(G \mid \gamma)$ is $\pi$-special.    

We now briefly review Isaacs' construction of ${\textup B}_{\pi}(G)$ in
\cite{bpi} and
Navarro's construction of the similar set, $\mbox{Irr}(G | Q, 1_Q)$.  In \cite{bpi} Isaacs
proves that if $G$ is $\pi$-separable, and if $\chi \in \mbox{Irr}(G)$, then
there is a unique (up to conjugacy) pair $(S, \varphi)$ maximal with
the property that $S \nrml \nrml G$, $\varphi \in \textup {Irr}(S)$ lies under $\chi$, and $\varphi$ is
$\pi$-factorable.  Such a pair $(S, \varphi)$ is called a maximal
factorable subnormal pair.  Denote by $T$ the stabilizer of $\varphi$ in the normalizer in $G$ of $S$.  Isaacs then shows that if $S \neq G$, then $T < G$
and induction defines a bijection from $\mbox{Irr}(T \mid \varphi)$ to $\mbox{Irr}(G \mid \varphi)$.  Therefore, we can let $\psi \in \textup{Irr}(T)$ be the unique character of $T$ lying over $\varphi$ such that $\psi^G = \chi$.  The
subnormal nucleus $(W, \gamma)$ of $\chi$ is then defined recursively
by defining $(W, \gamma)$ to be the subnormal nucleus of $\psi$.  (If $\chi$ is
$\pi$-factorable, then $(W, \gamma)$ is defined to be $(G, \chi)$.)  Note that by the
construction, $\gamma^G = \chi$ and $\gamma$ is $\pi$-factorable.
Also, it is shown that the subnormal nucleus of $\chi$ is
unique up to conjugacy.  The set $\textup{B}_{\pi}(G)$ is defined as the set of
irreducible characters of $G$ whose nucleus character $\gamma$
is $\pi$-special. 

In \cite{npiconst}, Navarro similarly defines the set $\mbox{Irr}(G | Q, 1_Q)$ for a $p$-solvable group $G$.  Navarro shows that if $\chi
\in \mbox{Irr}(G)$, then there is a unique pair $(N, \theta)$ maximal with
the property that $N \nrml G$, $\theta \in \textup{Irr}(N)$ lies under $\chi$, and $\theta$ is $p$-factorable.  Such a pair is called a maximal factorable normal
pair.  It is shown that if $N < G$, then $G_{\theta} < G$, and thus the Clifford correspondence implies that if $T = G_{\theta}$, there is a unique character $\psi \in \textup{Irr}(T \mid \theta)$ such that $\psi^G = \chi$.  Navarro then defines the normal nucleus $(U,
\epsilon)$ to be the normal nucleus of $(T, \psi)$.
(Again, if $\chi$ is $p$-factorable, then $(U, \epsilon) = (G,
\chi)$.)  Note that again $\epsilon^G = \chi$ and $\epsilon$ is
$p$-factorable.  If $Q$ is a Sylow $p$-subgroup of $U$, and if $\delta \in \textup{Irr}(Q)$ is defined by $\delta = (\epsilon_p)_Q$, then we say the pair $(Q, \delta)$ is a normal vertex for $\chi$, and it is shown that this
pair and the normal nucleus are unique up to conjugacy.  The set $\mbox{Irr}(G\mid Q,
1_Q)$ is defined as $\{\chi \in \irr(G)\mid \delta = 1_Q\}$, or
equivalently, the set of irreducible characters of $G$ with a
$p'$-special normal nucleus character.  Although Navarro only defines
the set $\irr(G | Q, \delta)$ when $\pi = p'$ and $G$ is $p$-solvable, the
same construction of the normal nucleus and vertex of a character
works if $\pi$ is an arbitrary set of primes and $G$ is
$\pi$-separable.  In this case, we will define the set $\textup{N}_{\pi}(G)$ to
be $\irr(G | Q, 1_Q)$, only now $Q$ is a Hall $\pi'$-subgroup
of the normal nucleus subgroup $U$ of $\chi$.  Thus ${\textup B}_{\pi}(G)$ consists of those irreducible characters of $G$ with a $\pi$-special subnormal nucleus character, and ${\textup N}_{\pi}(G)$ consists of those characters of $G$ with a $\pi$-special normal nucleus character.

Recall that if $\chi$ is any class function of the $\pi$-separable group $G$, then $\chi^0$ denotes the restriction of $\chi$ to the elements of $G$ whose order is a $\pi$-number.  Moreover, the set ${\textup I}_{\pi}(G)$ is a generalization of Brauer characters in a $p$-solvable group $G$ to a set of primes $\pi$ (so that if $\pi = p'$, and if $G$ is $p$-solvable, then the set ${\textup I}_{\pi}(G)$ is exactly the set of Brauer characters of $G$ for the prime $p$.).  In \cite{bpi} it is shown that if $\chi \in {\textup B}_{\pi}(G)$, then $\chi^0 \in {\textup I}_{\pi}(G)$, and in \cite{npiconst} it is shown that if $\eta \in {\textup N}_{\pi}(G)$, then $\eta^0 \in {\textup I}_{\pi}(G)$.  Both $\textup{B}_{\pi}(G)$ and $\textup{N}_{\pi}(G)$ are canonical
sets of lifts of $\textup{I}_{\pi}(G)$; in other words, if $\varphi \in
{\textup I}_{\pi}(G)$, then there is a unique character $\chi \in {\textup B}_{\pi}(G)$ and a unique
 character $\eta \in {\textup N}_{\pi}(G)$ such that $\chi^0 = \eta^0 = \varphi$.  Moreover,
Isaacs shows that if $\chi \in {\textup B}_{\pi}(G)$ and $N \nrml G$, then every
constituent $\theta$ of $\chi_N$ is in ${\textup B}_{\pi}(N)$.

We will need the following results.  In \cite{pispecial}, Isaacs
defines, for a subgroup $H$ of a $\pi$-separable group $G$ (where $2
\not \in \pi$), a linear
character $\delta_{(G:H)} \in \mbox{Irr}(H)$, called the standard sign character of
$H$.  This theorem lists some of the properties of $\delta_{(G:H)}$.

\medskip

{\bf {Theorem 1.1.}}  {\it Let $G$ be $\pi$-separable with $2 \not \in \pi$,
and let $\delta$ be the standard sign character for some subgroup $H$
of $G$.  Then the following hold:}

\medskip

\begin{itemize} 
\item[(1)]  {\it If $| G : H |$ is a $\pi$-number and $H$ is a maximal subgroup of
  $G$, then $\delta$ is the permutation sign character of the action
  of $H$ on the right cosets of $H$ in $G$.}
\item[(2)]  {\it $\mbox{core}_G(H) \subseteq {\textup{ker}}(\delta)$.}
\item[(3)]  {\it Suppose $H \sbs G$ has $\pi$ index.  If
$\psi \in \textup{Irr}(H)$ and $\psi^G = \chi \in \textup{Irr}(G)$, then $\chi$ is
$\pi$-special if and only if $\delta \psi$ is $\pi$-special.}

\end{itemize}

\medskip

Proof.   This is the content of Theorems 2.5 and B of
\cite{pispecial} and Lemma 2.1 of \cite{pipartial}.  $\Box$

\medskip

The aims of this paper, then, are threefold.  First, we prove the
statement made (without proof) in \cite{npiconst} that if $p$ is an odd
prime, the sets $\textup{B}_{p'}(G)$ and $\textup{N}_{p'}(G)$ coincide.  Secondly, we
give a counterexample to show that if $2 \in \pi'$, then $\textup{B}_{\pi}(G)$
and $\textup{N}_{\pi}(G)$ need not coincide.  Finally, we use the same
counterexample to show that if $\eta \in \textup{N}_{\pi}(G)$ and $M \nrml G$,
then it need not be the case that the constituents of $\eta_M$ are
themselves in $\textup{N}_{\pi}(M)$.

\medskip

{\bf{2.  Equality occurs if $G$ has odd order.}}  In this section we give a brief proof  that if
$| G |$ is odd or if $2 \in \pi$, then $\textup{B}_{\pi}(G) = \textup{N}_{\pi}(G)$.  This result was stated without proof in \cite{npiconst}.  Here
the field $\mathbb{Q}_{\pi}$ is defined by adjoining the $n$th roots of
unity to $\mathbb{Q}$ for every $\pi$-number $n$.

\medskip

{\bf {Lemma 2.1.}}  {\it Let $G$ be $\pi$-separable and assume $2 \in \pi$
or $| G |$ is odd.  Assume $\chi \in \textup{Irr}(G)$ has values in
$\mathbb{Q}_{\pi}$ and that $\chi^0 \in \ipig$.  Then $\chi \in \textup{B}_{\pi}(G)$.}

\medskip

Proof.  This is Theorem 12.3 of \cite{bpi}.  $\Box$
\medskip

{\bf {Lemma 2.2.}}  {\it Let $G$ be $\pi$-separable, and suppose $\eta \in
\textup{N}_{\pi}(G)$.  Then $\eta(g) \in \mathbb{Q}_{\pi}$ for all elements $g \in G$.}
\medskip

Proof.  Let $\sigma$ be any automorphism of $\mathbb{C}$ that
fixes the nth roots of unity for every $\pi$-number $n$.  Clearly
$\eta^{\sigma} \in \textup{N}_{\pi}(G)$ by the construction of $\textup{N}_{\pi}(G)$.
If $\eta^0 = \varphi \in \textup{I}_{\pi}(G)$, then since the values of $\varphi$ are
in $\mathbb{Q}_{\pi}$, then $\varphi^{\sigma} = \varphi$.  Since $\eta$ is the
unique lift of $\varphi$ in $\textup{N}_{\pi}(G)$, then it must be that
$\eta^{\sigma} = \eta$.  Thus the values of $\eta$ are in
$\mathbb{Q}_{\pi}$.  $\Box$

\medskip

{\bf {Corollary 2.3.}}  {\it If $G$ is $\pi$-separable and $2 \in \pi$ or
$|G| $ is odd, then $\textup{B}_{\pi}(G) = \textup{N}_{\pi}(G)$.}

\medskip

Proof.  Since $\textup{B}_{\pi}(G)$ and $\textup{N}_{\pi}(G)$ are both sets of
lifts of $\textup{I}_{\pi}(G)$, then $| \textup{B}_{\pi}(G)| = |\textup{N}_{\pi}(G)|$.  By Lemmas
2.1 and 2.2, $\textup{N}_{\pi}(G) \sbs \textup{B}_{\pi}(G)$.  Thus $\textup{N}_{\pi}(G) = \textup{B}_{\pi}(G)$.  $\Box$

\medskip

{\bf{3.  A counterexample of even order.}}  In this section we
construct the aforementioned counterexample to show that $\textup{B}_{\pi}(G)$ need not equal $\textup{N}_{\pi}(G)$, and we show that the constituents of the restriction of a character in $\npig$ to a normal subgroup $V$ need not be in $\npi(V)$.

Suppose $\Gamma$ is a finite group of order $n$, and let $p$ be a prime number.  Let $E$ be an elementary abelian group of order $p^n$.  Then we can let $\Gamma$ act on $E$ by associating to each element $x \in \Gamma$ one of the cyclic factors of $E$, and let the left multiplication action of $\Gamma$ on itself induce an action on the factors of $E$.  Let $G$ be the semidirect product of $\Gamma$ acting on $E$ with this action.  Then for each subgroup $L$ of $G$ such that $E \subseteq L \subseteq G$, we see that there is an irreducible character $\theta$ of $E$ such that $G_{\theta} = L$.  

\medskip

\textbf{Construction}:  Let $S_1$ be isomorphic to the
symmetric group on four elements, and let $A_1$ be isomorphic to the
alternating group on four elements.  Define the group $U_1$ as the
semidirect product of $S_1$ acting on $A_1$ with the conjugation
action, so that $U_1 = A_1 S_1$ and $S_1 \cap A_1 = 1$.  Let $K_1$
be the normal Klein four group in $A_1$, and note that $K_1 \triangleleft U_1$.  Define the
element $x \in S_1$ by setting $x$ equal to the permutation $(1 2)$.
Note that $\langle A_1, x \rangle \cong {\textup{Sym}}(4)$.  Let $H_1 = \langle
A_1, x \rangle$, and notice that $H_1$ is not subnormal in $U_1$,
$K_1 \triangleleft H_1$, and $H_1/K_1 \cong {\textup{Sym}}(3)$.
Let $L_1 \subseteq U_1$ be such that $K_1 \subseteq L_1 \subseteq H_1$
and $|H_1 : L_1| = 3$.  Note that $L_1 \cap A_1
= K_1$ and $L_1 A_1 = H_1$.  Moreover, note that ${\mathbb{O}}_3(U_1) = 1$.

Let $U_2$ be isomorphic to $U_1$, set $V_0 = U_1 \times U_2$, and define $\Gamma$ as the semidirect product
of ${\bf{C}}_2$ acting on $V_0$, with the nontrivial element of
${\bf{C}}_2$ acting by interchanging the components of $V_0$.  Note
that ${\textup{core}}_{\Gamma}(L_1) = 1$, and ${\mathbb{O}}_3(\Gamma) = 1$.

Let $|\Gamma| = n$, and let $E$ be an elementary abelian group of
order $3^n$.  Let $G$ be the semidirect product as defined in the paragraph preceding this construction, so that $G/E \cong \Gamma$.
Therefore there is a character $\theta \in \textup{Irr}(E)$ such that
$G_{\theta} = EL_1$.  Set $L = EL_1$, $A= EA_1$, $K = EK_1$, $H =
EH_1$, and $V = EV_0$.

Set $\pi = \{ 3 \}$ and note that $G$ is solvable.  Since $E$ is a
$3$-group, then $\theta$ is $\pi$-special.  Note that the $\pi$-special character $\theta$ extends to a unique
$\pi$-special character $\hat{\theta} \in \irr(K)$.  We see that since
$\theta$ and $\hat{\theta}$ uniquely determine each other, then
$A_{\hat{\theta}} = K$.
Therefore $\hat{\theta} \in \textup{Irr}(K)$ induces
irreducibly to a $\pi$-special character $\varphi \in \textup{Irr}(A)$.  Since
$\varphi$ is uniquely determined by $\theta$, and $L$ normalizes $A$,
then $L = G_{\theta} \subseteq G_{\varphi}$.

We now show that $G_{\varphi} = H$.  By a Frattini argument, we see that $G_{\varphi} \subseteq G_{\theta}A = H$.  Since $\varphi \in \textup{Irr}(A)$, then clearly $A \subseteq
G_{\varphi}$.  We showed in the previous paragraph that $L =
G_{\theta} \subseteq G_{\varphi}$.  Thus $H = AL \subseteq
G_{\varphi}$ and therefore $G_{\varphi} = H$.

We now claim that $(E, \theta)$ is a maximal factorable normal pair.  Recall that
${\mathbb{O}}_{3}(\Gamma)$ is trivial, and thus
${\mathbb{O}}_{3}(G/E)$ is trivial.  Suppose there exists a factorable normal pair $(N,
\psi)$ such that $E \subseteq N$ and $\theta$ lies under $\psi$, and suppose
$N/E$ is a nontrivial $2$-group.  Then the $\pi$-special factor of
$\psi$ must restrict to $\theta$, and thus $\theta$ must be invariant in
$N$.  Therefore $N \subseteq G_{\theta} = L$.  However,
${\textup{core}}_{G/E}(L/E)$ is trivial, and this yields a contradiction.
Thus $(E, \theta)$ is a maximal factorable normal pair.

We also claim that $(A, \varphi)$ is a maximal factorable subnormal pair.  Suppose
$(S, \sigma)$ is a factorable subnormal pair such that $A \triangleleft S$ and $\sigma$ lies over $\varphi$, and suppose $S/A$ is a nontrivial $2$-group.
Since $\varphi \in \irr(A)$ is $\pi$-special, $\varphi$ must be
invariant in $S$, and thus $S \subseteq G_{\varphi} = H$.  Since
$|H:A| = 2$ and $H$ is not subnormal in $G$, we have a contradiction.
Since there are no subnormal subgroups $T$ of $G$ such that $A
\triangleleft T$ and $T/A$ is a nontrivial $3$-group, then $(A,
\varphi)$ is a maximal factorable subnormal pair.

Note that since $|L:K| =2$ and the $\pi$-special character $\hat{\theta} \in \irr(K)$ is invariant in $L$, then $\hat{\theta}$
extends to a $\pi$-special character $\xi \in \irr(L)$, and we define $\eta \in
\irrg$ by $\eta = \xi^G$.  Note that $\theta$ lies under $\eta$, and
since $(E, \theta)$ is a maximal factorable normal pair, and the Clifford
correspondent $\xi \in \irr(G_{\theta}|\theta)$ for $\eta$ is
$\pi$-special, then $\eta \in \npig$.

We now show that $\eta$ is not in $\bpig$.  Note that
$\xi^H \in \irr(H)$, and since $\varphi \in \irr(A)$ is exactly $(\xi_K)^A = (\xi^H)_A$, then $\varphi$ lies under both $\xi^H$ and
$\eta$.  Recall that $(A,
\varphi)$ is a maximal factorable subnormal pair, and $H = G_{\varphi}$.  Since $|H:A|
= 2$ and $\varphi \in \irr(A)$ is $\pi$-special, every character of
$H$ lying over $\varphi$ must be $\pi$-factorable.  Thus $(H, \xi^H)$
is a subnormal nucleus for $\eta$, and $\eta \in \bpig$ if and
only if $\xi^H$ is $\pi$-special.

Let $\delta = \delta_{(H:L)} \in \irr(L)$ be the standard sign
character for $L \subseteq H$.  Note that $K \subseteq {\textup{ker}}(\delta)$.  Therefore we see that $\delta$ is the nonprincipal linear character of
$L/K$.  Since $\xi \in \irr(L)$ is $\pi$-special and $\delta \in
\irr(L)$ is $\pi'$-special, then $\xi \delta \in \irr(L)$ is not
$\pi$-special.  Thus $\xi^H$ is not
$\pi$-special by Theorem 1.1.  Therefore, since $(H, \xi^H)$ is a subnormal nucleus
for $\eta$, then $\eta$ is not in $\bpig$.

Finally, note that $(H, \xi^H)$ is a normal nucleus for $\xi^V$, which lies under $\eta$.  Since $\xi^H$ is not $\pi$-special, then $\xi^V$ is not $\pi$-special, and the constituents of $\eta_V$ are not in $\npi(V)$.  $\Box$

The above example raises some interesting questions.  Are there other families of characters which form lifts of the set $\ipig$?  Some results in this direction can be found in forthcoming papers of Mark Lewis.  What in general can be said about the set of lifts of a Brauer character of a $p$-solvable group $G$?  In a future paper, the author will describe some bounds on the number of lifts of a Brauer character of a $p$-solvable group.  There is still much to be known, though, about the set of lifts of a Brauer character of a $p$-solvable group. 
\medskip

{\bf {Acknowledgments}}

The above results are part of my doctoral thesis.  I am very grateful for the supportive environment of the University of Wisconsin-Madison Mathematics Department, the thoughtful suggestions of the referee, and especially the guidance of my advisor I. M. Isaacs.

\end{document}